\documentclass[10pt,draft,reqno]{amsart}
     \makeatletter
     \def\section{\@startsection{section}{1}%
     \z@{.7\linespacing\@plus\linespacing}{.5\linespacing}%
     {\bfseries
     \centering
     }}
     \def\@secnumfont{\bfseries}
     \makeatother
\newtheorem{theorem}{Theorem}[section]
\newtheorem{lemma}[theorem]{Lemma}
\newtheorem{proposition}[theorem]{Proposition}
\newtheorem{corollary}[theorem]{Corollary}
\theoremstyle{definition}
\newtheorem{definition}[theorem]{Definition}
\newtheorem{example}[theorem]{Example}
\theoremstyle{remark}
\newtheorem{remark}[theorem]{Remark}
\numberwithin{equation}{section}
\usepackage{amssymb,enumerate,bm}
\usepackage[numbers]{natbib}
\allowdisplaybreaks[4]

\newcommand{\ld}{\lambda}

\newcommand{\q}{\quad}

\newcommand{\eqd}{\overset{\mathrm d}{=}}

\newcommand{\wh}{\widehat}

\newcommand{\la}{\langle}
\newcommand{\ra}{\rangle}

\newcommand{\B}{\mathcal{B}}

\newcommand{\R}{\mathbb{R}}
\newcommand{\N}{\mathbb{N}}
\newcommand{\Z}{\mathbb{Z}}

\newcommand{\rd}{{\mathbb R^d}}
\newcommand{\law}{\mathcal L}
\newcommand{\sek}{\int_0^{\infty}}

\newcommand{\plim}{\mathop{\mbox{\rm p-$\lim$}}}
\def\1{\bm{1}}

\begin{document}

\title[Semi-selfdecomposable distributions and related processes]
{Stochastic integral characterizations of semi-selfdecomposable distributions and\\
related Ornstein-Uhlenbeck type processes}

\author{Makoto Maejima}
\address{Makoto Maejima:
Department of Mathematics, Keio University, 3-14-1, Hiyoshi, Kohoku-ku, Yokohama 223-8522, Japan}
\email{maejima@math.keio.ac.jp}

\author{Yohei Ueda}
\address{Yohei Ueda:
Department of Mathematics, Keio University, 3-14-1, Hiyoshi, Kohoku-ku, Yokohama 223-8522, Japan}
\email{ueda@2008.jukuin.keio.ac.jp}



\begin{abstract}
In this paper, three topics on semi-selfdecomposable distributions are studied.
The first one is to characterize semi-selfdecomposable distributions by stochastic integrals 
with respect to L\'evy processes.
This characterization defines a mapping from an infinitely divisible
distribution with finite $\log$-moment to 
a semi-selfdecomposable distribution.
The second one is to introduce and study a Langevin type equation and the corresponding 
Ornstein-Uhlenbeck
type process whose limiting distribution is semi-selfdecomposable.
Also, semi-stationary Ornstein-Uhlenbeck type processes with semi-selfdecomposable 
distributions are constructed.
The third one is to study the iteration of the mapping above.
The iterated mapping is expressed as a single mapping with a different integrand. 
Also, nested subclasses of 
the class of semi-selfdecomposable distributions are considered, and
it is shown that the limit of these nested subclasses is the closure of the class of 
semi-stable distributions. 
\end{abstract}

\maketitle

\section{Introduction}
Let $I(\rd)$ be the class of all infinitely divisible distributions on $\rd$ and let $\{X_t^{(\mu)},t\geq 0\}$ be an 
$\rd$-valued L\'evy process
with $\mu\in I(\rd)$ as its distribution at time 1.
Many subclasses of $I(\rd)$ have recently been investigated in many aspects.
Among those, there are characterizations of those classes in terms of stochastic integrals with
respect to L\'evy processes.
In such cases, we define mappings
$$
\Phi_f(\mu) = \law \left ( \sek f(t)dX_t^{(\mu)}\right ) ,\quad \mu\in\mathfrak D(\Phi_f)\subset I(\rd)
$$
for nonrandom measurable functions $f\colon [0,\infty)\rightarrow \R$,
where $\law (X)$ is the law of a random variable $X$ and $\mathfrak D (\Phi_f)$ is the domain of 
a mapping $\Phi_f$ that is
the class of $\mu\in I(\rd)$ for which $\sek f(t)dX_t^{(\mu)}$ is definable.
For the definition of stochastic integrals with respect to L\'evy processes of nonrandom measurable 
functions, see the next section. 
When we consider the composition of two mappings $\Phi_f$ and $\Phi_g$, denoted by $\Phi_g\circ\Phi_f$,
the domain of $\Phi_g\circ\Phi_f$ is $\mathfrak D (\Phi_g\circ \Phi_f) = 
\{ \mu\in I(\rd) : \mu\in \mathfrak D (\Phi_f)\,\,
\text{and}\,\, \Phi_f(\mu)\in \mathfrak D (\Phi_g)\}$.
Once we define such a mapping, we can characterize a subclass of $I(\rd)$ as the range of $\Phi_f$, 
$\mathfrak R (\Phi_f)$, say.
Among such classes, there are the Jurek class, the class of selfdecomposable distributions, 
the Goldie-Steutel-Bondesson class,
the Thorin class, the class of generalized type $G$ distributions and so on.
(For details on these, see, e.g., \citet{MaejimaSato2009}.)
Also, by iterating a mapping $\Phi_f$, we can define a sequence of nested subclasses 
$\mathfrak R (\Phi_f^m), m\in\N$,
where $\Phi_f^m$ is the $m$ times composition of the mapping $\Phi_f$ itself.

The class of selfdecomposable distributions, denoted by $L(\rd)$, has the longest history 
in the study of subclasses of $I(\rd)$.
Let $\wh \mu (z), z\in\rd$, be the characteristic function of $\mu$.
$\mu\in I(\rd)$ is said to be selfdecomposable if for any $b>1$, there exists a distribution 
$\rho_b$ such that $\wh \mu (z) = \wh\mu (b^{-1}z)\wh{\rho}_b (z)$.
This $\rho_b$ automatically belongs to $I(\rd)$.
$\mu\in L(\rd)$ is also a limiting distribution of normalized partial sums of independent 
random variables under infinitesimal condition, and has the stochastic integral representation 
with respect to a L\'evy process, 
which is $\law\left (\sek e^{-t}dX_t\right)$ with $\law(X_1)\in I_{\log}(\rd)$,
where $I_{\log}(\rd) = \{ \mu\in I(\rd) \colon \int _{\rd} \log ^+|x| \mu (dx) <\infty \}$, 
$\log^+|x|= (\log|x|)\vee 0$, and $|x|$ is the Euclidean norm of $x\in\rd$.
Furthermore, $\mu\in L(\rd)$ is the limiting distribution of 
the solution of a Langevin equation with L\'evy noise.
More precisely, let $\{X_t, t\ge 0\}$ be a L\'evy process on $\rd$, $c\in\R$, 
and let $M$ be an $\rd$-valued random variable.
The Langevin equation is
\begin{equation}\label{lang}
Z_t = M + X_t -c \int_0^t Z_sds, \quad t\geq 0,
\end{equation}
and the following is known, (see, e.g., \citet{Sato's_book2003}).
$$
Z_t = e^{-ct}M + e^{-ct}\int_0^te^{cs}dX_s, \quad t\ge 0,
$$
is an almost surely unique solution of \eqref{lang},
and if $c>0$, $E\left[\log ^+|X_1|\right ] <\infty$, and $M$ is independent of $\{X_t\}$,
then $\law(Z_t)\to \mu \in L(\rd)$ as $t\rightarrow \infty $.
We also know that, for a fixed $c>0$, the equation
\begin{equation}\label{lang_stationary}
Z_t-Z_s=X_t-X_s-c\int_s^t Z_udu, \quad -\infty<s\leq t<\infty,
\end{equation}
has an almost surely unique stationary solution 
$$
Z_t=\int_{-\infty}^te^{-c(t-u)}dX_u,\quad t\in\R,
$$
where $\{X_t,t\in\R\}$ is a h.-i.s.r.m.-process (whose precise definition is given in Section 2) satisfying $E\left[\log ^+|X_1|\right ] <\infty$.
(See, e.g., \citet{Sato's_book2003} or \citet{MaejimaSato2003}.)
This stationary solution fulfills that $\law(Z_t)=\law\left(\int_0^\infty e^{-cu}dX_u\right)
\in L(\rd)$ for all $t\in\R$.
Also, it is recognized that some selfdecomposable distributions on $\R$ are very important 
in the area of mathematical finance,
(see \citet{Yor2007}).

In \citet{MaejimaNaito1998}, the concept of the selfdecomposability was extended to the
semi-selfdecomposability.
Here, $\mu\in I(\rd)$ is called {\it semi-selfdecomposable} if there exist $b>1$ 
and $\rho\in I(\rd)$ such that $\wh \mu (z) = \wh\mu (b^{-1}z)\wh\rho(z)$.
We call this $b$ a span of $\mu$, and we denote the class of 
all semi-selfdecomposable distributions with span $b$ by $L(b, \rd)$.
From the definitions, $L(b,\rd)\supsetneqq L(\rd)$ and $L(\rd) = \bigcap _{b>1}L(b,\rd)$.
$\mu\in L(b,\rd)$ is also realized as a limiting distribution of normalized partial sums of 
independent random variables under infinitesimal condition when the limit is taken through 
a geometric subsequence.
A typical example is a semi-stable distribution, where $\mu\in I(\rd)$ is said to be 
semi-stable with span $b$ if 
there exist $a>1$ and $c\in\rd$ satisfying $\wh\mu(z)^a=\wh\mu(bz)e^{i\la c,z\ra}$.
Recently, several natural examples of semi-selfdecomposable distributions have appeared in the literature.
We will mention some of them in the next section.

In \citet{MaejimaSato2003}, they gave a stochastic integral characterization of 
$\mu\in L(b,\rd)$ in terms of,
not L\'evy process, but natural semi-L\'evy process.
Here a semi-L\'evy process with period $p>0$ is an additive process with periodically 
stationary increments with period $p$ and
natural additive process was defined in \citet{Sato2004} as semimartingale additive process in terms of
the L\'evy-Khintchine triplet.
Namely, they showed that for each $b>1$, $\mu\in L(b,\rd)$ if and only if 
$\mu=\law\left( \sek e^{-t}dX_t\right )$, where $\{X_t\}$ is a semi-L\'evy process with
period $p=\log b$ and $\law(X_p)\in I_{\log}(\rd)$.
Our first topic of this paper is to give a stochastic integral characterization of 
$\mu\in L(b,\rd)$ in terms of L\'evy process.
If all natural semi-L\'evy processes can be expressed as stochastic integrals with respect to
L\'evy processes, this problem is trivial from a result in \citet{MaejimaSato2003} just mentioned now.
However, as we will see in Example \ref{example} later, it is not the case.
Once we could solve this problem, we would define a mapping $\Phi_b$ from $\mathfrak D (\Phi_b)$
into $I(\rd)$ and we can enjoy many stories similar to those about $L(\rd)$.
For instance, we can characterize $L(b,\rd)$ as the range of the mapping $\Phi_b$, that is,
$L(b, \rd) = \Phi_b(I_{\log}(\rd))$.

Our second topic  is to construct and study a Langevin type equation and
the corresponding Ornstein-Uhlenbeck type processes related to semi-selfdecomposable distributions,
which are analogies of \eqref{lang} and \eqref{lang_stationary} in the case of selfdecomposable
distributions, not in terms of semi-L\'evy processes given in \citet{MaejimaSato2003}, 
but in terms of L\'evy processes.
Namely, we introduce a Langevin type equation and give its unique solution, 
which we call an Ornstein-Uhlenbeck type process.
We then show that the limit of the Ornstein-Uhlenbeck type process exists in law,
when the noise process has finite $\log$-moment, and the limiting distribution is 
semi-selfdecomposable.
We also construct semi-stationary Ornstein-Uhlenbeck type process whose 
marginal distributions are semi-selfdecomposable.

Our third topic is to look for the ranges of the iterated mappings $\Phi_b^m$ and its limit.
In \citet{MaejimaNaito1998}, the nested subclasses of $L(b,\rd)$, $L_m(b,\rd), m\in\Z_+$, 
are defined as follows:
$\mu \in L_m(b,\rd)$ if and only if there exists  $\rho\in L_{m-1}(b,\rd)$ such that
$\wh \mu (z) = \wh\mu (b^{-1}z)\wh\rho (z)$, where $L_0(b, \rd) = L(b,\rd)$.
Then we will show
\begin{equation}\label{nested_subclass_of_ssd}
L_m(b, \rd) = \Phi^{m+1} _b\left(I_{\log ^{m+1}}(\rd)\right),\quad m\in \Z_+,
\end{equation}
where $I_{\log^{m+1}}(\rd) = \{ \mu\in I(\rd) \colon \int _{\rd} 
(\log ^+|x|)^{m+1} \mu (dx) <\infty \}$.
The relation \eqref{nested_subclass_of_ssd} implies
that the limit of these nested subclasses is the closure of the class of 
semi-stable distributions,
where the closure is taken under convolution and weak convergence.

Organization of this paper is the following.
In Section 2, we explain some notation and give preliminaries and 
some examples of semi-selfdecomposable distributions.
In Section 3, the first topic is considered.
In Sections 4--6, we study the second topic.
Finally, in Section 7, we treat the third topic.
\section{Notation, preliminaries and examples}
In this section, we explain necessary notation, and give some preliminaries and examples.

Let $J$ be $\R$ or $[0,\infty)$, and $\B_J^0$ the class of all bounded Borel sets in $J$.
An $\rd$-valued independently scattered random measure (abbreviated as i.s.r.m.) $X=\{X(B),B\in\B_J^0\}$
is said to be homogeneous if $\law(X(B))=\law(X(B+a))$ for all $B\in\B_J^0$ and $a\in\R$ satisfying $B+a\in\B_J^0$.
See \citet{MaejimaSato2003} and \citet{Sato2004,Sato2006a},
for the definition and deep study of stochastic integrals
of nonrandom measurable functions $f\colon J\to\R$
with respect to $\rd$-valued i.s.r.m.'s $X$, 
denoted by $\int_B f(s)X(ds), B\in\B_J^0$.
For a fixed $t_0\in J$, we use the symbol
$$
\int_{t_0}^t f(s)X(ds)=
\begin{cases}
\int_{(t_0,t]} f(s)X(ds),&\text{for }t\in (t_0,\infty),\\
0,&\text{for }t=t_0,\\
-\int_{(t,t_0]} f(s)X(ds),&\text{for }t\in J\cap(-\infty,t_0),
\end{cases}
$$
which is understood to be a c\`adl\`ag modification,
(see Remark 3.16 of \citet{MaejimaSato2003}).
If $\{X_t,t\geq 0\}$ is a L\'evy process, then there exists 
a unique $\rd$-valued homogeneous i.s.r.m.\,\,$X$ over $[0,\infty)$ satisfying $X_t=X([0,t])$ a.s.\,\,for each $t\geq 0$.
Then $\int_0^tf(s)dX_s$ is defined by 
$\int_0^tf(s)X(ds)$ for $t\in[0,\infty)$. 
The improper stochastic integral $\sek f(s)dX_s$ is defined as the limit in probability of 
$\int_0^tf(s)dX_s$
as $t\to\infty$ whenever the limit exists.
See also \citet{Sato2006b}.
In this paper, we say that a stochastic process $\{X_t,t\in\R\}$ on $\rd$ is a \emph{h.-i.s.r.m.-process}
if there exists an $\rd$-valued homogeneous i.s.r.m.\,\,$X$ over $\R$ such that
$X_t=\int_0^tX(du),t\in\R$.
We define a stochastic integral $\int_s^tf(u)dX_u,-\infty<s\leq t<\infty$ 
of a nonrandom measurable function $f\colon\R\rightarrow \R$ 
with respect to this process $\{X_t, t\in\R\}$
by $\int_s^tf(u)X(du)$.
See also \citet{Sato's_book2003} and \citet{MaejimaSato2003}.
The improper stochastic integral $\int_{-\infty}^tf(u)dX_u$ is defined as the limit in probability of 
$\int_s^tf(u)dX_u$ as $s\rightarrow -\infty$, provided that this limit exists.
If the improper stochastic integral $\int_{-\infty}^tf(u)dX_u$
is definable for $t\in\R$, then we regard it as a c\`adl\`ag process, 
since such a modification always exists.

Throughout this paper, we use the L\'evy-Khintchine representation of the characteristic function of $\mu\in I(\R^d)$ 
in the following way:
$$
\widehat{\mu}(z)=\exp\left\{-\frac 12 \la z,Az\ra +i\la \gamma,z\ra +\int_\rd 
\left(e^{i\la z,x\ra}-1-\frac{i\la z,x\ra}{1+|x|^2}\right)\nu(dx)\right\},\quad z\in\rd,
$$
where $\la\cdot,\cdot\ra$ denotes Euclidean inner product on $\rd$, $A$ is a nonnegative-definite symmetric $d\times d$ 
matrix, $\gamma\in\rd$, 
and $\nu$ is a measure, called L\'evy measure, satisfying $\nu(\{0\})=0$ and $\int_\rd(|x|^2\wedge 1)\nu(dx)<\infty$. 
We call $(A,\nu,\gamma)$ the L\'evy-Khintchine triplet of $\mu$ and we write $\mu=\mu_{(A,\nu,\gamma)}$ 
when we want to emphasize the L\'evy-Khintchine triplet.
$C_{\mu}(z),z\in\rd$, denotes the cumulant function of $\mu\in I(\rd)$, that is, $C_{\mu}(z)$ is the unique continuous function
satisfying $\wh\mu(z) = e^{C_{\mu}(z)}$ and $C_{\mu}(0)=0$.
When a random variable $X$ has the distribution $\mu$, we sometime write $C_X(z)$ for $C_{\mu}(z)$.

We also use the polar decomposition \eqref{polar} of the L\'evy measure $\nu$ of $\mu\in I(\rd)$ with $0<\nu(\R^d)\le\infty$.
There exist a measure $\ld$ on $S:=\{x\in\rd\colon|x|=1\}$ with $0<\ld(S)\le\infty$ and
a family $\{\nu_{\xi}, \xi\in S\}$ of measures on $(0,\infty)$ such that
$\nu_{\xi}(B)$ is measurable in $\xi$ for each $B\in\mathcal B((0,\infty))$,
$0<\nu_{\xi}((0,\infty))\le\infty$ for each $\xi\in S$ and
\begin{align}\label{polar}
\nu(B)=\int_S \ld(d\xi)\int_0^{\infty} \1_B(r\xi)\nu_{\xi}(dr),\quad
B\in \mathcal B (\mathbb R^d \setminus \{ 0\}).
\end{align}
Here $\ld$ and $\{\nu_{\xi}\}$ are uniquely determined by $\nu$
up to multiplication of measurable functions $c(\xi)$ and $c(\xi)^{-1}$, respectively, with $0<c(\xi)<\infty$.
We say that $\nu$ has the polar decomposition $(\ld ,\nu_{\xi})$, and
$\ld$ and $\nu_{\xi}$ are called the spherical and the radial components of $\nu$, respectively.
(See, e.g., \citet{Barndorff-NielsenMaejimaSato2006}, Lemma 2.1.)

Recently, several natural examples of semi-selfdecomposable distributions have appeared in the literature.
In \citet{Watanabe2002}, he showed that the distribution of a certain supercritical branching process 
and the first hitting time of Brownian motion starting at the origin on the unbounded Sierpinski gasket
on $\mathbb R^2$ are both semi-selfdecomposable.
Also, let $\{N_t,t\ge 0\}$ be a Poisson process and $\{X_t\}$ a L\'evy process on $\rd$ independent of $\{N_t\}$. 
Suppose $E\left [\log^+|X_1|\right ] <\infty$ and $b>1$.
Then $\law \left (\sek b^{-N_{t-}}dX_t\right )$ is semi-selfdecomposable with span $b$.
(Theorem 3.2 of \citet{KondoMaejimaSato2006}.)
In a recent paper by \citet{LindnerSato2009}, we can also find several examples of semi-selfdecomposable 
distributions with the form $\law\left (\sek  c^{-N_{t-}}dX_t\right )$, where $\{(N_t, X_t)\}$ is a bivariate
compound Poisson process with L\'evy measure concentrated on the three points $(1,0), (0,1)$ and $(1,1)$.
Another recent example is found by \citet{Pacheco-Gonzalez2009} in some financial modeling.
These indicate introducing of semi-selfdecomposable distributions allows us more flexibility in
stochastic modeling.
\section{Stochastic integral characterizations of semi-selfdecomposable distributions}

As we mentioned in Introduction, in this section, we introduce a mapping from a subset of $I(\rd)$ into $I(\rd)$, 
by which semi-selfdecomposable distributions can be characterized.
\begin{definition}
Let $b>1$ and $\mu\in I(\R^d)$. Define a mapping $\Phi_b$ by
\begin{equation}\label{def_of_Phi_b}
\Phi_b(\mu):= \law\left(\int_0^\infty b^{-[t]}dX_t^{(\mu)}\right),
\end{equation}
provided that this improper stochastic integral is definable, where $[x]$ denotes the largest integer not 
greater than $x\in\R$.
\end{definition}
The domain of the mapping $\Phi_b$, where the improper stochastic integral in \eqref{def_of_Phi_b} is definable, 
is given as follows by Theorem 2.4 of \citet{Sato2006}.
\begin{proposition}\label{domain_of_Phib}
$\mathfrak{D}(\Phi_b)=I_{\log}(\rd)$.
\end{proposition}
We start with the following theorem. 
\begin{theorem}\label{equivalence_between_decomposability_and_mapping}
Fix any $b>1$. Let $\mu$ and $\rho$ be distributions on $\rd$. Then,
\begin{equation}\label{relation_c.f.}
\rho\in I(\rd)\quad\text{and}\quad \widehat{\mu}(z)=\widehat{\mu}(b^{-1}z)\widehat{\rho}(z)
\end{equation}
if and only if
\begin{equation}\label{mapping}
\rho\in I_{\log}(\rd)\quad\text{and}\quad \mu=\Phi_b(\rho).
\end{equation}
\end{theorem}
\begin{proof}
To show the ``if" part, suppose \eqref{mapping}. 
Note that
\begin{equation}\label{product}
\Phi_b(\rho) = \law \left (\sek b^{-[t]}dX_t^{(\rho)}\right ) = \law\left (\sum _{j=0}^{\infty}b^{-j} \left(X_{j+1}^{(\rho)} - X_{j}^{(\rho)}\right)\right ).
\end{equation}
Then
$$
\wh\mu(z)  = \prod_{j=0}^{\infty}\wh\rho (b^{-j}z) = \prod_{j=1}^{\infty}\wh\rho (b^{-j}z)\times \wh\rho (z)
 = \prod_{k=0}^{\infty}\wh\rho (b^{-k}(b^{-1}z))\times \wh\rho (z) = \wh\mu (b^{-1}z)\wh\rho(z),
$$
which concludes \eqref{relation_c.f.}.

We next show the ``only if" part. Assume \eqref{relation_c.f.}. 
Then as can be seen in \citet{Wolfe1983}, we have
$$
\wh\mu(z) = \wh\mu (b^{-1}z)\wh\rho(z) = \wh\mu (b^{-2}z)\wh\rho(b^{-1}z)\wh\rho(z)
 = \cdots  = \wh\mu(b^{-n}z)\prod_{j=0}^{n-1}\wh\rho (b^{-j}z),
$$
for all $n\in\N$.
Hence it follows that $\prod_{j=0}^{\infty}\wh\rho (b^{-j}z)$ exists and equals $\wh\mu(z)$,
which implies $\rho\in I_{\log}(\rd)$ by \citet{Wolfe1983}.
Then $\Phi_b(\rho)$ is definable and 
satisfies \eqref{product}.
Thus we have $\Phi_b(\rho)=\mu$, which yields \eqref{mapping}.
\end{proof}
Theorem \ref{equivalence_between_decomposability_and_mapping} yields the following.
\begin{corollary}\label{s.s.d.}
Fix any $b>1$. Then, the range $\mathfrak R (\Phi_b)$ is the class of all semi-selfdecomposable distributions 
with span $b$ on $\rd$, namely,
$$
\Phi_b\left(I_{\log}(\rd)\right)=L(b,\R^d).
$$
\end{corollary}
The injectivity of the mapping $\Phi_b$ is shown as follows.
\begin{proposition}\label{ssd_injective} For each $b>1$, 
the mapping $\Phi_b$ is injective.
\end{proposition}
\begin{proof}
Let $\rho_1,\rho_2\in I_{\log}(\rd)$ and $\mu=\Phi_b(\rho_1)=\Phi_b(\rho_2)$.
Then, Theorem \ref{equivalence_between_decomposability_and_mapping} yields that
$$
\widehat{\mu}(z)=\widehat{\mu}(b^{-1}z)\widehat{\rho}_1(z)=\widehat{\mu}(b^{-1}z)\widehat{\rho}_2(z).
$$
Since $\widehat{\mu}(b^{-1}z)\neq 0$ for all $z\in\rd$ by the infinite divisibility of $\mu$, it follows that 
$\widehat{\rho}_1(z)=\widehat{\rho}_2(z)$.
\end{proof}
\begin{remark}
As mentioned in Introduction, if $\mu\in L(\rd)$, then there exists $\mu_0\in I_{\log}(\rd)$ such that 
$\mu=\law\left ( \sek e^{-t}dX_t^{(\mu_0)}\right )$, and it is known that this $\mu_0$ is uniquely determined by $\mu$.
We have just shown that if $\mu\in L(b, \rd)$, then there exists $\mu_b\in I_{\log}(\rd)$ such that 
$\mu=\law\left ( \sek b^{-[t]}dX_t^{(\mu_b)}\right )$, and the uniqueness of $\mu_b$ is assured by 
Proposition \ref{ssd_injective}.
If $\mu\in L(\rd)$, then $\mu\in L(b,\rd)$ for any $b>1$.
Then it is natural to ask what relation there is between $\mu_0$ and $\mu_b$ with $b>1$.
We answer this question below.
Fix $b>1$. Let $\mu=\mu_{(A,\nu,\gamma)}$, $\mu_0={\mu_0}_{(A_0,\nu_0,\gamma_0)}$ and $\mu_b={\mu_b}_{(A_b,\nu_b,\gamma_b)}$. Denote the polar decompositions of $\nu$, $\nu_0$ and $\nu_b$ by $(\lambda, \nu_{\xi})$, $(\lambda_0,\nu_{0,\xi})$ and $(\lambda_b,\nu_{b,\xi})$, respectively.
Note that $\mu\in L(\rd)$ if and only if 
$$
\nu_{\xi}(dr)=\frac{k_\xi(r)}{r}dr,\quad r>0,
$$
where $k_\xi(r)$ is a nonnegative function, which is measurable in $\xi$, and is nonincreasing 
and right-continuous in $r$. 
(See \citet{Sato's_book1999}, Theorem 15.10.)
We have
\begin{align}
\label{A_0}A&=\int_0^\infty e^{-2t}A_0dt=2^{-1}A_0,\\
\notag\gamma&=\int_0^\infty e^{-t}dt\left\{\gamma_0+\int_\rd x\left(\frac{1}{1+e^{-2t}|x|^2}-\frac{1}{1+|x|^2}\right)
\nu_0(dx)\right\}\\
\label{gamma_0}&=\gamma_0+\int_{\rd\setminus\{0\}}x\left(\frac{\arctan |x|}{|x|}-\frac{1}{1+|x|^2}\right)\nu_0(dx),
\end{align}
and it follows from Theorem 41 (ii) of \citet{Sato's_book2003} that $\lambda_0=\lambda$ and $\nu_{0,\xi}(dr)=-dk_\xi(r)$, 
up to multiplication of positive finite measurable functions $c(\xi)$ and $c(\xi)^{-1}$.
On the other hand, Theorem \ref{equivalence_between_decomposability_and_mapping} yields that $\mu_b$ is 
an infinitely divisible distribution satisfying 
$\widehat\mu(z)=\widehat\mu(b^{-1}z)\widehat\mu_b(z)$. Therefore
\begin{align}
\label{A_b}A_b&=\left(1-b^{-2}\right)A,\\
\label{gamma_b}\gamma_b&=\left(1-b^{-1}\right)\gamma-\int_{\rd}x\left(\frac{1}{1+|x|^2}-\frac{1}{1+|bx|^2}\right)\nu(b\,dx),
\end{align}
and
\begin{align*}
\nu_b(B)&=\nu(B)-\nu(bB)=\int_S\lambda(d\xi)\int_0^\infty\1_B(r\xi)\frac{k_\xi(r)-k_\xi(br)}{r}dr\\
&=\int_S\lambda(d\xi)\int_0^\infty\1_B(r\xi)\frac{\nu_{0,\xi}\left((r,br]\right)}{r}dr,\qquad B\in\mathcal B_0(\rd).
\end{align*}
Then, it follows that
\begin{equation}\label{nu_0_and_nu_b}
\lambda_b=\lambda_0=\lambda\quad\text{and}\quad\nu_{b,\xi}(dr)=\frac{\nu_{0,\xi}\left((r,br]\right)}{r}dr=
\frac{k_\xi(r)-k_\xi(br)}{r}dr\ \lambda\text{-a.e.}\,\,\xi\in S,
\end{equation}
up to multiplication of positive finite measurable functions $c(\xi)$ and $c(\xi)^{-1}$.
One can see the relation between $\mu_0$ and $\mu_b$ by \eqref{A_0}, \eqref{gamma_0}, \eqref{A_b}, \eqref{gamma_b} 
and \eqref{nu_0_and_nu_b}.
\end{remark}
As also mentioned in Introduction, \citet{MaejimaSato2003} characterized semi-selfdecomposable distributions 
by stochastic integrals with respect to natural semi-L\'evy processes.
The following theorem is another version of Corollary \ref{s.s.d.} in this paper and  Corollary 5.4 of 
\citet{MaejimaSato2003}, and connects them. 
For $b>1$, let $G_b$ denote the totality of bounded periodic measurable functions with period $\log b$.
\begin{theorem}\label{period}
Fix any $b>1$. Then we have
$$
L(b,\R^d)=\left\{\law\left(\int_0^\infty e^{-t}g(t)dX_t^{(\mu)}\right)\colon g\in G_b \text{ and } 
\mu\in I_{\log}(\rd)
\right\}.
$$
\end{theorem}
\begin{proof}
Let $\widetilde\mu\in L(b,\rd)$. Then Corollary \ref{s.s.d.} yields that $\widetilde\mu=\Phi_b(\mu)$ 
for some $\mu\in I_{\log}(\rd)$. If we let
$$
g(t):=b^{\frac{t}{\log b}-\left[\frac{t}{\log b}\right]},
$$
then $g\in G_b$ and $e^{-t}g(t)=b^{-[t/\log b]}$. 
It follows that
$$
\widetilde\mu=\Phi_b(\mu)
=\law\left(\int_0^\infty b^{-[t/\log b]}dX_{t/\log b}^{(\mu)}\right)
=\law\left(\int_0^\infty e^{-t}g(t)dX_t^{\left(\mu^{1/\log b}\right)}\right),
$$
where for $p>0$, $\mu^p$ is an infinitely divisible distribution with characteristic function $\widehat{\mu}(z)^p$.

Conversely, suppose that $\widetilde\mu=\law\left(\int_0^\infty e^{-t}g(t)dX_t^{\left(\mu\right)}\right)$ 
for $g\in G_b$ and $\mu\in I_{\log}(\rd)$.
Putting $Y_t:=\int_0^tg(s)dX_s^{(\mu)}$, we have $\widetilde\mu=\law\left(\int_0^\infty e^{-t}dY_t\right)$ 
due to Theorem 4.6 of \citet{Sato2004}, and we see that $\{Y_t\}$ is a natural semi-L\'evy process with 
period $\log b$.
Moreover, we have $E\left[\log^+|Y_{\log b}|\right]<\infty$ since the L\'evy measure of $\law(Y_{\log b})$ 
denoted by $\widetilde \nu$ satisfies that
\begin{align*}
\int_{|x|>1}\log|x|\widetilde\nu(dx)&=\int_0^{\log b}ds\int_\rd\log^+|g(s)x|\nu(dx)\\
&\leq \log b\int_\rd\log^+\biggl(\sup_{s\in[0,\log b]}|g(s)||x|\biggr)\nu(dx)<\infty,
\end{align*}
where $\nu$ is the L\'evy measure of $\mu$. Then Corollary 5.4 of \citet{MaejimaSato2003} implies that 
$\widetilde\mu\in L(b,\rd)$.
\end{proof}
In the proof of Theorem \ref{period}, $Y_t = \int_0^tg(s)dX_s^{(\mu)}$ with $g\in G_b$ and 
$\mu\in I_{\log}(\rd)$ is proved to be a natural semi-L\'evy process.
However, any natural semi-L\'evy process is not necessarily expressed in this way as is shown 
in the following example.
\begin{example}\label{example}
We claim that not all natural semi-L\'evy processes can be expressed as stochastic integrals with respect 
to L\'evy processes.
Fix an arbitrary $p>0$. Let $\varphi\colon [0,p]\rightarrow [0,\infty)$ be a nondecreasing function satisfying 
$\varphi(0)=0$ which is continuous but not absolutely continuous, for example, Cantor's function on $[0,p]$.
Suppose that $\mu\in I(\rd)\setminus\{\delta_0\}$.
Then there exists a semi-L\'evy process $\{Y_t,t\geq 0\}$ with period $p$ such that 
$C_{Y_t}(z)={\varphi(t)}C_{\mu}(z)$ for $z\in\rd$ and $t\in [0,p]$, due to 
Proposition 2.2 of \citet{MaejimaSato2003}.
Furthermore, it follows from the monotonicity of $\varphi$ that $\{Y_t\}$ is natural.
However, we cannot express $\{Y_t\}$ in the form that $Y_t=\int_0^tf(s)dX_s$ for any measurable function $f$ and any 
L\'evy process $\{X_t\}$.
Indeed, if $\int_0^tf(s)dX_s$ is definable for all $t\in[0,\infty)$, then 
$C_{\int_0^tf(s)dX_s}(z)=\int_0^tC_{X_1}(f(s)z)ds$ which is absolutely continuous in $t$, 
although $C_{Y_t}(z)$ is not absolutely continuous in $t$ by the property of $\varphi$.
\end{example}
\section{A Langevin type equation}

The purpose of this and the following two sections is to find a Langevin type equation like \eqref{lang} 
or \eqref{lang_stationary} related to semi-selfdecomposable distributions.
The ideas of proofs of the results below come from Sections 2.2 and 2.4 of \citet{Sato's_book2003} 
and \citet{MaejimaSato2003}. 

For this purpose, we first consider the following \emph{Langevin type equation}:
\begin{equation}\label{Langevin_t_0}
Z_t = M + X_{[ct]/c}-X_{[ct_0]/c} -(b-1) \int_{t_0}^{t} Z_sd[cs], \quad t\geq t_0,
\end{equation}
where $t_0\in\R$, $c>0$, $b>1$, $\{X_t, t\in\R\}$ is a h.-i.s.r.m.-process on $\rd$,
 $M$ is an $\rd$-valued random variable, 
and $\int_0^{t} Z_sd[cs]$ is as follows: $d[cs]$ denotes the Lebesgue-Stieltjes measure associated with $s\mapsto [cs]$, 
which is equal to $\sum_{k\in\Z}\delta_{k/c}(ds)$, and $\int_{(\alpha,\beta]} f(s)d[cs]$, 
written as $\int_\alpha^\beta f(s)d[cs]$, exists for $-\infty<\alpha<\beta<\infty$ and any 
measurable function $f$, since it is a finite sum in fact. 
A stochastic process $\{Z_t\}$ is said to be a \emph{solution} of the Langevin type equation \eqref{Langevin_t_0} 
or an \emph{Ornstein-Uhlenbeck type process} generated by $\{X_t\}$, $b$ and $c$ starting from $Z_{t_0}=M$ 
if sample paths of $\{Z_t\}$ are right-continuous with left limits and $\{Z_t\}$ satisfies \eqref{Langevin_t_0}  almost surely.
We claim that
\begin{equation}\label{solution_t_0}
Z_t=b^{-([ct]-[ct_0])}M+b^{-[ct]}\int_{[ct_0]/c}^{[ct]/c}b^{[cs]}dX_s, \q t\ge t_0,
\end{equation}
is an almost surely unique solution of \eqref{Langevin_t_0}.
\begin{theorem}\label{thm_Langevin_t_0}
Suppose that $\{X_t, t\in\R\}$ is a h.-i.s.r.m.-process on $\rd$, $t_0\in\R$, $c>0$, $b>1$, and $M$ is 
an $\rd$-valued random variable. 
Then, $\{Z_t\}$ in \eqref{solution_t_0} is an almost surely unique solution of the equation \eqref{Langevin_t_0}.
\end{theorem}
\begin{proof} 
If we define $\{Z_t\}$ by \eqref{solution_t_0}, then it is a c\`adl\`ag process. 
If $t_0\leq t<t_0+1/c$, then $Z_t=M$ and it satisfies \eqref{Langevin_t_0}. 
Let $t\geq t_0+1/c$. Then, it follows that
\begin{align*}
(b&-1)\int_{t_0}^{t} Z_sd[cs]\\
&=(b-1)\sum_{k=[ct_0]+1}^{[ct]}Z_{k/c}
=(b-1)\sum_{k=[ct_0]+1}^{[ct]}\left\{b^{-(k-[ct_0])}M+b^{-k}\int_{[ct_0]/c}^{k/c}b^{[cs]}dX_s\right\}\\
&=\left(1-b^{-([ct]-[ct_0])}\right)M +(b-1)\sum_{k=[ct_0]+1}^{[ct]}b^{-k}\sum_{\ell=[ct_0]+1}^kb^{\ell-1}
\left(X_{{\ell}/{c}}-X_{{(\ell-1)}/{c}}\right)\\
&=\left(1-b^{-([ct]-[ct_0])}\right)M +(b-1)\sum_{\ell=[ct_0]+1}^{[ct]}b^{\ell-1}\left(X_{{\ell}/{c}}
-X_{({\ell-1)}/{c}}\right)
\sum_{k=\ell}^{[ct]}b^{-k}\\
&=\left(1-b^{-([ct]-[ct_0])}\right)M +\sum_{\ell=[ct_0]+1}^{[ct]}\left(1-b^{-[ct]+\ell-1}\right)\left(X_{{\ell}/{c}}
-X_{{(\ell-1)}/{c}}\right)\\
&=\left(1-b^{-([ct]-[ct_0])}\right)M +X_{[ct]/c}-X_{[ct_0]/c}-b^{-[ct]}\int_{[ct_0]/c}^{[ct]/c}b^{[cs]}dX_s\\
&=M+X_{[ct]/c}-X_{[ct_0]/c}-Z_t.
\end{align*}
This yields \eqref{Langevin_t_0}.

It remains to prove the uniqueness of the solution of \eqref{Langevin_t_0}. 
Suppose that two $\{Z^{(1)}_t\}$ and $\{Z^{(2)}_t\}$ are the solutions of \eqref{Langevin_t_0}. 
Setting $Z_t:=Z^{(1)}_t-Z^{(2)}_t$, we have
\begin{equation}\label{uniqueness}
Z_t=-(b-1)\int_{t_0}^tZ_sd[cs]\qquad\text{for } t\geq t_0,\quad\text{a.s.}
\end{equation}
Let us show that
\begin{equation}\label{uniqueness_induction}
Z_t=0,\qquad \text{for }\,\,\left (t_0\vee \frac{[ct_0]+n-1}{c}\right )\leq t<\frac{[ct_0]+n}{c},\quad\text{a.s.}
\end{equation}
for any $n\in\N$ by induction. \eqref{uniqueness_induction} is true for $n=1$, since the right-hand side 
of \eqref{uniqueness} is zero for $t_0\leq t<([ct_0]+1)/c$. 
Assume that \eqref{uniqueness_induction} holds for $n=1,2,\dots,m$. 
Then, for $([ct_0]+m)/c\leq t<([ct_0]+m+1)/c$, \eqref{uniqueness} can be reduced to that
$$
Z_t=-(b-1)\int_{([ct_0]+m-1)/c}^{([ct_0]+m)/c}Z_sd[cs]=-(b-1)Z_{([ct_0]+m)/c},
$$
which implies that $X_{([ct_0]+m)/c}=0$ by letting $t=([ct_0]+m)/c$ and thus $Z_t=0$ for 
$([ct_0]+m)/c\leq t<([ct_0]+m+1)/c$. 
Hence \eqref{uniqueness_induction} is true for $n=m+1$. 
Therefore it holds with probability one that for any $t\geq 0$, $Z_t=Z^{(1)}_t-Z^{(2)}_t=0$.
\end{proof}
\section{Limiting distributions of Ornstein-Uhlenbeck type processes}
In this section, we study the Langevin type equation \eqref{Langevin_t_0} with $t_0=0$:
\begin{equation}\label{Langevin}
Z_t = M + X_{[ct]/c} -(b-1) \int_0^{t} Z_sd[cs], \quad t\geq 0,
\end{equation}
where $\{X_t,t\geq 0\}$ is a L\'evy process on $\rd$, $c>0$, $b>1$, and $M$ is an $\rd$-valued random variable.
Theorem \ref{thm_Langevin_t_0} yields that
\begin{equation}\label{solution}
Z_t=b^{-[ct]}M+b^{-[ct]}\int_0^{[ct]/c}b^{[cs]}dX_s, \quad t\geq 0,
\end{equation}
is an almost surely unique solution of \eqref{Langevin}. 
We show that if $M$ is independent of $\{X_t\}$ and $\law (X_1)\in I_{\log}(\rd)$, then the limit of 
$\law(Z_t)$ exists 
as $t\rightarrow \infty $ and the limiting distribution is semi-selfdecomposable with span $b$. 
\begin{theorem}\label{limiting_distribution}
Suppose that $\{X_t, t\geq 0\}$ is a L\'evy process on $\rd$, $c>0$, $b>1$, and $M$ is an $\rd$-valued random 
variable independent of $\{X_t\}$. 
Then, $\law(X_1)\in I_{\log}(\rd)$ if and only if $Z_t$ in \eqref{solution} converges in law to some 
$\rd$-valued random variable
as $t\rightarrow \infty $. 
This limit $\lim_{t\rightarrow \infty }\law(Z_t)$ is equal to $\law\left(\int_0^\infty b^{-[cs]-1}dX_s\right)$ 
which is semi-selfdecomposable with span $b$ and does not depend on $M$. 
Furthermore, if we let $\law(M):=\lim_{t\rightarrow \infty }\law(Z_t)$, then $\law(Z_t)=\law(M)$ for all $t\geq 0$.
\end{theorem}
\begin{proof}
In the rest of the paper, we write $\wh{\law}(X)(z) $ for the characteristic function of $\law(X)$ for notational simplicity.
It follows that
\begin{align}
\notag \widehat{\law}(Z_t)(z)
&=\widehat{\law}(M)\left(b^{-[ct]}z\right)\exp\left\{\int_0^{[ct]/c}C_{X_1}
\left(b^{\left[cs-[ct]\right]}z\right)ds\right\}\\
\notag&=\widehat{\law}(M)\left(b^{-[ct]}z\right) \exp\left\{\int_0^{[ct]/c}
C_{X_1}\left(b^{[-cu]}z\right)du\right\}\\
\label{sol_cf}&=\widehat{\law}(M)\left(b^{-[ct]}z\right) \exp\left\{\int_0^{[ct]/c}
C_{X_1}\left(b^{-[cu]-1}z\right)du\right\}.
\end{align}
Note that $\widehat{\law}(M)\left(b^{-[ct]}z\right)$ tends to $1$ as $t\rightarrow \infty $.

Assume $\law(X_1)\in I_{\log}(\rd)$. 
Then, \eqref{sol_cf} tends to the cumulant function of the improper stochastic 
integral $\int_0^\infty b^{-[cu]-1}dX_u$, whose law is semi-selfdecomposable with span $b$ due 
to Corollary \ref{s.s.d.}. 
Also, $\lim_{t\rightarrow \infty }\law(Z_t)=\law\left(\int_0^\infty b^{-[cu]-1}dX_u\right)$ does not 
depend on $M$ and if we set $\law(M):=\law\left(\int_0^\infty b^{-[cu]-1}dX_u\right)$, then \eqref{sol_cf} is
\begin{align*}
C_{Z_t}(z) 
& =\int_0^{\infty}C_{X_1}\left(b^{-\left[cu+[ct]\right]-1}z\right)du
+\int_0^{[ct]/c}C_{X_1}\left(b^{-[cu]-1}z\right)du\\
&=\int_{[ct]/c}^{\infty}C_{X_1}\left(b^{-[cv]-1}z\right)dv
+\int_0^{[ct]/c}C_{X_1}\left(b^{-[cu]-1}z\right)du\\
&=\int_0^{\infty}C_{X_1}\left(b^{-[cu]-1}z\right)du,
\end{align*}
which yields that $\law(Z_t)=\law\left(\int_0^\infty b^{-[cu]-1}dX_u\right)$ for all $t\geq 0$.

Next assume that $Z_t$ converges in law as $t\rightarrow \infty $, namely, $\int_0^{n/c}C_{X_1}
\left(b^{-[cu]-1}z\right)du$ tends to the cumulant function of some $\mu\in I(\rd)$ with L\'evy measure 
$\nu$ as $n\rightarrow \infty $. 
Denoting the L\'evy measure of $\law(X_1)$ by $\nu_{X_1}$, we have
$$
\lim_{k\rightarrow \infty}\int_0^{n_k/c}du\int_\rd\left(|b^{-[cu]-1}x|^2\wedge 1\right)\nu_{X_1}(dx)=
\int_\rd(|x|^2\wedge 1)\nu(dx)
$$
for some subsequence, due to the proof of Theorem 8.7 in \citet{Sato's_book1999}. It follows from the monotone 
convergence theorem that
$\int_0^{\infty}du\int_\rd\left(|b^{-[cu]-1}x|^2\wedge 1\right)\nu_{X_1}(dx)<\infty$, which implies 
$\law(X_1)\in I_{\log}(\rd)$ by virtue of Lemma 2.7 of \citet{Sato2006}.
\end{proof}
Theorem \ref{limiting_distribution} yields that if $\law(X_1)\in I_{\log}(\rd)$, then $Z_t$ in \eqref{solution} 
converges in law as $t\rightarrow \infty $. 
Then, it might be natural to ask whether or not $Z_t$ converges in probability as $t\rightarrow \infty $.
The following proposition is the answer.
\begin{proposition}
Suppose that $\{X_t, t\geq 0\}$ is a L\'evy process on $\rd$, $c>0$, $b>1$, and $M$ is an $\rd$-valued random 
variable independent of $\{X_t\}$. If $\law(X_1)$ is not any $\delta$-distribution, then $Z_t$ in 
\eqref{solution} does 
not converge in probability as $t\rightarrow \infty $.
\end{proposition}
\begin{proof}
Suppose that $\mu:=\law(X_1)$ is not any $\delta$-distribution.
Then, Lemma 13.9 of \citet{Sato's_book1999} yields that $\left|\widehat{\mu}(z_0)\right|<1$ for some $z_0\in\rd$.
As
$$
Z_t-Z_{t-1/c}=b^{-[ct]}\int_{[ct-1]/c}^{[ct]/c}b^{[cs]}dX_s+(b^{-[ct]}-b^{-[ct-1]})
\left(M+\int_0^{[ct-1]/c}b^{[cs]}dX_s\right)
$$
by \eqref{solution}, we have
\begin{align*}
&\left|\widehat{\law} \left(Z_t-Z_{t-1/c}\right)(z)\right|\leq \left|\widehat{\law}\left(b^{-[ct]}
\int_{[ct-1]/c}^{[ct]/c}b^{[cs]}dX_s\right)(z)\right|\\
&\qquad\qquad=\left|\widehat{\law}\left(b^{-1}\left(X_{[ct]/c}-X_{[ct-1]/c}\right)\right)(z)\right|
=\left|\widehat{\mu}(b^{-1}z)\right|^{1/c}
\end{align*}
for any $t\geq 1/c$ and $z\in\rd$.
This yields that for all $t\ge 1/c$,
$$
\left|\widehat{\law} \left(Z_t-Z_{t-1/c}\right)(bz_0)\right|\leq 
\left|\widehat{\mu}(z_0)\right|^{1/c}<1.
$$
Then $Z_t-Z_{t-1/c}$ does not tends to zero in probability as $t\rightarrow \infty$.
Thus $Z_t$ does not converge in probability as $t\rightarrow \infty $.
\end{proof}
The following remark is about the relation between the Langevin type equation \eqref{Langevin} and the mapping 
$\Phi_b$.
\begin{remark}\label{relation_between_Langevin_and_mapping}
Fix $b>1$ and $c>0$. Let $\mu\in I_{\log}(\rd)$. Consider the limiting solution in law 
$\lim_{t\rightarrow \infty }\law(Z_t)$ 
of the Langevin type equation \eqref{Langevin} with a L\'evy process $\{X_t\}$ satisfying $\law(b^{-1}X_{1/c})=\mu$ and an $\rd$-valued random variable $M$ independent of $\{X_t\}$. 
It follows from Theorem \ref{limiting_distribution} that
\begin{align*}
\Phi_b(\mu)&
=\law\left(\int_0^\infty b^{-[t]}dX_t^{(\mu)}\right)
=\law\left(\int_0^\infty b^{-[t]}dX_t^{\left(\law (b^{-1}X_{1/c})\right)}\right)\\
&=\law\left(\int_0^\infty b^{-[t]}d(b^{-1}X_{t/c})\right)
=\law\left(\int_0^\infty b^{-[cs]-1}dX_s\right)=\lim_{t\rightarrow \infty }\law(Z_t).
\end{align*}
Thus the mapping $\Phi_b$ can be defined also as the limiting distribution of the solution of 
the Langevin type equation 
\eqref{Langevin}.
\end{remark}
We conclude this section with the Markov property of our Ornstein-Uhlenbeck type processes.
\begin{proposition}
Suppose that $\{X_t, t\geq 0\}$ is a L\'evy process on $\rd$, $c>0$, $b>1$, and $M$ is an $\rd$-valued random 
variable independent of $\{X_t\}$. 
Then, the process $\{Z_t\}$ in \eqref{solution} is a Markov process satisfying
\begin{equation}\label{Markov}
P\left(Z_t\in B\mid Z_s=x\right)=P\left(b^{-([ct]-[cs])}x+\int_0^{([ct]-[cs])/c}b^{-[cu]-1}dX_u\in B\right)
\end{equation}
for $0\leq s\leq t$ and $B\in\B(\rd)$.
\end{proposition}
\begin{proof}
Since
$$
Z_t=b^{-([ct]-[cs])}Z_s+b^{-[ct]}\int_{[cs]/c}^{[ct]/c}b^{[cu]}dX_u,
$$
we can easily see the Markov property of $\{Z_t\}$ by virtue of the independent increment property of the L\'evy process 
$\{X_t\}$.
\eqref{Markov} is shown as follows:
\begin{align*}
E&\left[e^{i\la z,Z_t\ra}\bigm| Z_s=x\right]\\
&=\exp\left\{i\la z,b^{-([ct]-[cs])}x\ra+\int_{[cs]/c}^{[ct]/c}C_{X_1}\left (b^{[cu-[ct]]}z\right )du\right\}\\
&=\exp\left\{i\la z,b^{-([ct]-[cs])}x\ra+\int_{0}^{([ct]-[cs])/c}C_{X_1}\left (b^{[-cv]}z\right)du\right\}\\
&=\exp\left\{i\la z,b^{-([ct]-[cs])}x\ra+\int_{0}^{([ct]-[cs])/c}C_{X_1}\left (b^{-[cv]-1}z\right )du\right\}.
\end{align*}
\end{proof}
\section{Semi-stationary Ornstein-Uhlenbeck type processes having semi-selfdecomposable distributions}
This section is concerned with the following \emph{Langevin type equation} which has similar properties to 
\eqref{lang_stationary}:
\begin{equation}\label{Langevin_semi-stationary}
Z_t-Z_s =X_{[ct]/c}-X_{[cs]/c} -(b-1) \int_s^{t} Z_ud[cu], \quad -\infty<s\leq t<\infty,
\end{equation}
where $c>0$, $b>1$, and $\{X_t,t\in\R\}$ is a h.-i.s.r.m.-process on $\rd$.
A stochastic process $\{Z_t\}$ is said to be a \emph{solution} of the Langevin type equation 
\eqref{Langevin_semi-stationary} 
or an \emph{Ornstein-Uhlenbeck type process} generated by $\{X_t\}$, $b$ and $c$ if sample paths of 
$\{Z_t\}$ are right-continuous with left limits and $\{Z_t\}$ satisfies \eqref{Langevin_semi-stationary} 
almost surely.

We show that if $E[\log^+|X_1|]<\infty$, then $\int_{-\infty}^tb^{[cu]}dX_u$ is definable for each $t\in\R$ and
\begin{equation}\label{solution_semi-stationary}
Z_t=b^{-[ct]}\int_{-\infty}^{[ct]/c}b^{[cu]}dX_u, \quad t\in\R,
\end{equation}
is an almost surely unique semi-stationary solution of \eqref{Langevin_semi-stationary},
where the \emph{semi-stationarity} 
of $\{Z_t\}$ means $\{Z_{t+p}\}\eqd\{Z_t\}$ for a fixed $p>0$.
Here $\eqd$ stands for equality in all finite-dimensional distributions.
This $p$ is called the 
\emph{period} of the semi-stationary 
process $\{Z_t\}$.

To prove this, we prepare two lemmas.
\begin{lemma}\label{lemma_semi-stationary_Langevin}
Let $\{X_t,t\in\R\}$ be a h.-i.s.r.m.-process on $\rd$. Suppose that $b>1$ and $c>0$.
Then, the following three statements are equivalent:
\begin{enumerate}[\rm (i)]
\item $\law(X_1)\in I_{\log}(\rd)$,
\item $\int_{-\infty}^0b^{[ct]}dX_t$ is definable,
\item there exists an Ornstein-Uhlenbeck type process $\{Z_t\}$ generated by $\{X_t\}$, $b$ and $c$ satisfying 
$\plim_{t\rightarrow -\infty}b^{[ct]}Z_t=0$, where $\plim$ stands for limit in probability.
\end{enumerate}
If {\rm (iii)} holds, then $\{Z_t\}$ with the properties in {\rm (iii)} is almost surely unique and expressed as 
\eqref{solution_semi-stationary}.
\end{lemma}
\begin{proof}
We first show that (i) and (ii) are equivalent. 
Theorem 2.4 of \citet{Sato2006} yields that $\law(X_1)=\law\left(-X_{(-1)-}\right)\in I_{\log}(\rd)$ if and only if 
$\int_0^\infty b^{[-cu]}d\left(-X_{(-u)-}\right)$ is definable.
Lemma 4.8 of \citet{MaejimaSato2003} implies that $\int_0^\infty b^{[-cu]}d\left(-X_{(-u)-}\right)$ is definable 
if and only if $\int_{-\infty}^0b^{[cu]}dX_u$ is definable.
Thus (i) and (ii) are equivalent.

We next show that (ii) implies (iii). Assume that (ii) holds. Then, $\{Z_t\}$ in \eqref{solution_semi-stationary} 
is definable. 
It satisfies \eqref{Langevin_semi-stationary} due to Theorem \ref{thm_Langevin_t_0} by letting $t_0=s$ and 
$Z_{t_0}=M=b^{-[cs]}\int_{-\infty}^{[cs]/c}b^{[cu]}dX_u$. It follows from \eqref{solution_semi-stationary} 
and (ii) that
$$
\plim_{t\rightarrow -\infty}b^{[ct]}Z_t=\plim_{t\rightarrow -\infty}\int_{-\infty}^{[ct]/c}b^{[cu]}dX_u=0.
$$

Finally, we show that (iii) implies (i), the uniqueness of $\{Z_t\}$ in (iii), and the expression \eqref{solution_semi-stationary}.
If (iii) holds, then $\{Z_t\}$ in (iii) satisfies 
\eqref{Langevin_semi-stationary}. 
Theorem \ref{thm_Langevin_t_0} yields that for each $s\in \R$, with probability one,
$$
Z_t=b^{-([ct]-[cs])}Z_s+b^{-[ct]}\int_{[cs]/c}^{[ct]/c}b^{[cu]}dX_u,\qquad\text{for } t\ge s,
$$
namely, with probability one,
$$
b^{[ct]}Z_t-b^{[cs]}Z_s=\int_{[cs]/c}^{[ct]/c}b^{[cu]}dX_u,\qquad\text{for } t\ge s.
$$
By letting $s\rightarrow -\infty$, it follows from (iii) that for each $t\in\R$,
$$
b^{[ct]}Z_t=\plim_{k\rightarrow \infty}\int_{-k/c}^{[ct]/c}b^{[cu]}dX_u,\qquad\text{a.s.}
$$
By a similar argument to that in the proof of Theorem \ref{limiting_distribution}, we have $\law(X_1)\in I_{\log}(\rd)$. 
Thus (i) holds. Then (ii) holds and it follows that for each $t\in\R$,
$$
b^{[ct]}Z_t=\int_{-\infty}^{[ct]/c}b^{[cu]}dX_u,\qquad\text{a.s.}
$$
However, since the both sides of the equation above have c\`adl\`ag paths, we have, with probability one,
$$
b^{[ct]}Z_t=\int_{-\infty}^{[ct]/c}b^{[cu]}dX_u,\qquad\text{for all }t\in\R.
$$
This yields the almost sure uniqueness of $\{Z_t\}$ in (iii), and the expression \eqref{solution_semi-stationary}.
\end{proof}
\begin{lemma}\label{semi-stationarity_iff_mildness}
Let $\{X_t,t\in\R\}$ be a h.-i.s.r.m.-process on $\rd$, $b>1$ and $c>0$. 
Suppose that $\{Z_t\}$ is an 
Ornstein-Uhlenbeck type process generated by $\{X_t\}$, $b$ and $c$. 
Then, $\{Z_t\}$ is semi-stationary 
if and only if $\plim_{t\rightarrow -\infty}b^{[ct]}Z_t=0$. Let these conditions be fulfilled. 
Then, semi-stationary process 
$\{Z_t\}$ has a period $1/c$. Moreover, $\{Z_t\}$ with the properties above is almost surely unique and expressed 
as \eqref{solution_semi-stationary}.
\end{lemma}
\begin{proof}
We first show the ``only if" part.
Suppose that $\{Z_t\}$ is a semi-stationary Ornstein-Uhlenbeck type process generated by $\{X_t\}$, $b$ and $c$ 
with period $p>0$.
Since $\{Z_t\}$ has c\`adl\`ag paths, for any sequence $\{t_n,n\in\N\}\subset[0,p]$, 
there exists its subsequence 
$\{t_{n_k},k\in\N\}$ satisfying $Z_{t_{n_k}}$ converges almost surely to some $\rd$-valued random variable 
as $k\to\infty$. 
This implies the relative compactness of $\{\law(Z_t)\colon t\in[0,p]\}$ which is equal to 
$\{\law(Z_t)\colon t\in\R\}$ 
by the semi-stationarity of $\{Z_t\}$. Hence $\{\law(Z_t)\colon t\in\R\}$ is tight by Prohorov's theorem. 
Then, it follows that for any $\varepsilon>0$, 
$$
P\left(|b^{[ct]}Z_t|>\varepsilon\right)\leq \sup_{s\in\R}P\left(|Z_s|>b^{-[ct]}\varepsilon\right)\to 0,\quad\text{as }t\to -\infty.
$$
Thus $\plim_{t\rightarrow -\infty}b^{[ct]}Z_t=0$.

We next show the ``if" part.
Suppose that $\{Z_t\}$ is an Ornstein-Uhlenbeck type process generated by $\{X_t\}$, $b$ and $c$ 
satisfying $\plim_{t\rightarrow -\infty}b^{[ct]}Z_t=0$. Then $\{Z_t\}$ has the form 
\eqref{solution_semi-stationary} due to 
Lemma \ref{lemma_semi-stationary_Langevin}. Let $-\infty<t_1<t_2<\dots<t_n<\infty$. Then, for each $j=2,3,\dots,n$, we have
\begin{align*}
b^{[ct_j]}Z_{t_j+1/c}-b^{[ct_{j-1}]}Z_{t_{j-1}+1/c}
&=b^{-1}\int_{[ct_{j-1}+1]/c}^{[ct_j+1]/c}b^{[cu]}dX_u\\
&=\int_{[ct_{j-1}]/c}^{[ct_j]/c}b^{[cv]}d(X_{v+1/c}-X_{1/c}),
\end{align*}
which is equal in law to
$$
\int_{[ct_{j-1}]/c}^{[ct_j]/c}b^{[cv]}dX_v
=b^{[ct_j]}Z_{t_j}-b^{[ct_{j-1}]}Z_{t_{j-1}}.
$$
Since $\{b^{[ct]}Z_{t+1/c}\}$ and $\{b^{[ct]}Z_t\}$ have independent increment property due to the expression 
\eqref{solution_semi-stationary}, it follows that $\{b^{[ct]}Z_{t+1/c}\}\eqd \{b^{[ct]}Z_t\}$, which yields $\{Z_{t+1/c}\}\eqd \{Z_t\}$.
This is the semi-stationarity of $\{Z_t\}$ with period $1/c$.

The almost sure uniqueness of $\{Z_t\}$ follows from Lemma \ref{lemma_semi-stationary_Langevin}.
\end{proof}
Now, we prove the following theorem on the relation between semi-stationary Ornstein-Uhlenbeck type processes 
and semi-selfdecomposable distributions.
\begin{theorem}\label{thm_semi-stationary_Langevin}
Suppose that $\{X_t, t\in\R\}$ is a h.-i.s.r.m.-process on $\rd$, $c>0$, and $b>1$. 
Then, $\law(X_1)\in I_{\log}(\rd)$ if and only if \eqref{Langevin_semi-stationary} has a semi-stationary solution. 
If $\law(X_1)\in I_{\log}(\rd)$, $\{Z_t\}$ in \eqref{solution_semi-stationary} is an almost surely unique 
semi-stationary solution of \eqref{Langevin_semi-stationary} and it satisfies 
$\law(Z_t)=\law\left(\int_0^\infty b^{-[cu]-1}dX_u\right)\in L(b,\rd)$ for all $t\in \R$. 
In this case, the semi-stationary process $\{Z_t\}$ has a period $1/c$.
\end{theorem}
\begin{proof}
The most parts of this theorem have already been proved in Lemmas \ref{lemma_semi-stationary_Langevin} and 
\ref{semi-stationarity_iff_mildness}.
It remains to show the statement that if $\law(X_1)\in I_{\log}(\rd)$, then $\law(Z_t)
=\law\left(\int_0^\infty b^{-[cu]-1}dX_u\right)\in L(b,\rd)$ for all $t\in \R$. 
However, it follows that for all $t\in\R$,
\begin{align*}
C_{Z_t}(z)&=\int_{-\infty}^{[ct]/c}C_{X_1}\left(b^{\left[cu-[ct]\right]}z\right)du
=\int_0^\infty C_{X_1}\left(b^{[-cv]}z\right)dv\\
&=\int_0^\infty C_{X_1}\left(b^{-[cv]-1}z\right)dv
=C_{\int_0^\infty b^{-[cv]-1}dX_v}(z),
\end{align*}
which is the cumulant function of some semi-selfdecomposable distribution with span $b$, due to Corollary \ref{s.s.d.}. 
\end{proof}
The following remark, 
which is similar to Remark \ref{relation_between_Langevin_and_mapping},
is about the relation between the Langevin type equation \eqref{Langevin_semi-stationary} and 
the mapping $\Phi_b$.
\begin{remark}
Fix $b>1$ and $c>0$. Let $\mu\in I_{\log}(\rd)$. Consider the almost surely unique semi-stationary Ornstein-Uhlenbeck 
type process $\{Z_t\}$ generated by a h.-i.s.r.m.-process $\{X_t\}$ satisfying $\law(b^{-1}X_{1/c})=\mu$, and $b$ and $c$.
Theorem \ref{thm_semi-stationary_Langevin} and the same calculations as those in Remark 
\ref{relation_between_Langevin_and_mapping} yield that
$$
\Phi_b(\mu)=\law(Z_t),\quad\text{for all }t\in\R.
$$
Hence the mapping $\Phi_b$ can be defined also as the distribution of the semi-stationary solution of 
the Langevin type 
equation \eqref{Langevin_semi-stationary}.
\end{remark}
\section{Nested subclasses of $L(b,\rd)$ given by iterating the mapping $\Phi_b$}
We now go back to the mapping $\Phi_b$ itself again.
The iterated mapping of $\Phi_b$ can be expressed by one stochastic integral as follows.
\begin{theorem}\label{iteration}
Suppose $m\in\Z_+$. The domain of $\Phi_b^{m+1}$ is
\begin{equation}\label{iteration domain}
\mathfrak{D}(\Phi_b^{m+1})=I_{\log^{m+1}}(\R^d).
\end{equation}
Let
$$
f_m(u):= \int_0^u\binom{[v]+m}{m}dv,
$$
and let $f_m^*$ be its inverse function. Then
$$
\Phi_b^{m+1}(\mu)=\law\left(\int_0^\infty b^{-[f_m^*(t)]}dX_t^{(\mu)}\right),\quad \text{for }
\mu\in I_{\log^{m+1}}(\R^d).
$$
\end{theorem}
\begin{remark}
If we let
$$
\widetilde{f}_m(u):= \int_0^u\frac{v^m}{m!}dv,
$$
then its inverse function is $\widetilde{f}_m^*(t)=\{(m+1)!\,t\}^{\frac{1}{m+1}}$ and $e^{-\widetilde{f}_m^*(t)}$ 
is the integrand of the stochastic integral of the iteration of the mapping in the case of $L(\R^d)$, 
(see Remark 58 of \citet{Sato's_book2003}).
One can see the difference between $f_m$ and $\widetilde{f}_m$ by
$$f_m(u)= \int_0^u\binom{[v]+m}{m}dv=\int_0^u\frac{([v]+1)([v]+2)\cdots ([v]+m)}{m!}dv.$$
\end{remark}
\vskip 3mm
\begin{proof}[Proof of Theorem \ref{iteration}]
We prove the statement by induction. If $m=0$, the assertion is true by the definition of $\Phi_b$ and 
Proposition \ref{domain_of_Phib}. Assume that the assertion is true for $0,1,\dots,m-1$ in place of $m$. 
Let
$$
\widetilde{\Phi}_{b,m+1}(\mu):=\law\left(\int_0^\infty b^{-[f_m^*(t)]}dX_t^{(\mu)}\right).
$$
Then $\mathfrak{D}^0(\widetilde{\Phi}_{b,m+1})=\mathfrak{D}(\widetilde{\Phi}_{b,m+1})=I_{\log^{m+1}}(\R^d)$ 
due to Proposition 4.3 of \citet{Sato2006b}, where $\mathfrak{D}^0(\Phi_f)$ denotes the set of all $\mu\in I(\rd)$
 satisfying $\int_0^\infty|C_{\mu}(f(t)z)|dt<\infty$. 
If $\mu\in I_{\log^{m+1}}(\R^d)\subset I_{\log^{m}}(\R^d)$, 
then $\Phi_b^{m}(\mu)=\law\left(\int_0^\infty b^{-[f_{m-1}^*(t)]}dX_t^{(\mu)}\right)$ 
by the assumption of induction, and
\begin{align}
\notag
\int_0^\infty&\left|C_{\Phi_b^{m}(\mu)}(b^{-[t]}z)\right|dt\\
\notag&\leq  \int_0^\infty dt\int_0^\infty \left|C_{\mu}\left(b^{-[f_{m-1}^*(s)]-[t]}z\right)\right|ds\\
\notag&=\int_0^\infty dt\int_0^\infty \left|C_{\mu}\left(b^{-\left[f_{m-1}^*(s)+[t]\right]}z\right)\right|ds\\
\notag&=\int_0^\infty dt\int_{[t]}^\infty \left|C_{\mu}(b^{-[u]}z)
\right|\binom{[u]{-}[t]{+}m{-}1}{m{-}1}du\\
\notag&=\int_0^\infty\left|C_{\mu}(b^{-[u]}z)\right|du\int_{0}^{[u]+1}
\binom{[u]{-}[t]{+}m{-}1}{m{-}1}dt\\
\notag&=\int_0^\infty\left|C_{\mu}(b^{-[u]}z)\right|\binom{[u]+m}{m}du\\
\label{Phi^m}&=\int_0^\infty\left|C_{\mu}(b^{-[f_m^*(t)]}z)\right|dt<\infty,
\end{align}
since $\mu\in I_{\log^{m+1}}(\R^d)=\mathfrak{D}^0(\widetilde{\Phi}_{b,m+1})$. 
Note that we have used above the formula
$$\int_0^{n-k+1}\binom{n-[t]}{k}dt
=\sum_{j=0}^{n-k}\binom{n-j}{k}
=\binom{n+1}{k+1}\quad\text{for }k\leq n.$$
Hence $I_{\log^{m+1}}(\R^d)\subset \mathfrak{D}(\Phi_b^{m+1})$. 
By similar calculations to \eqref{Phi^m}, we have
$$
\int_0^\infty C_{\Phi_b^{m}(\mu)}(b^{-[t]}z)dt=\int_0^\infty C_{\mu}(b^{-[f_m^*(t)]}z)dt,
\quad\text{for }\mu\in I_{\log^{m+1}}(\rd),
$$
where the use of Fubini's theorem is permitted by the finiteness of \eqref{Phi^m}. Thus
$$\Phi_b^{m+1}(\mu)=\widetilde{\Phi}_{b,m+1}(\mu),\quad\text{for }\mu\in I_{\log^{m+1}}(\rd).$$
To conclude \eqref{iteration domain}, it remains to prove that $\mathfrak{D}(\Phi_b^{m+1})\subset I_{\log^{m+1}}(\rd)$. If $\mu\in I(\rd)$ 
with L\'evy measure $\nu$ satisfies $\mu\notin I_{\log^{m+1}}(\rd)$, there exists $n\in\{0,1,\dots,m\}$ 
such that $\mu\in I_{\log^{n}}(\rd)\setminus I_{\log^{n+1}}(\rd)$ (consider $I_{\log^{0}}(\rd)$ to be $I(\rd)$). 
If $n=0$, $\mu\notin I_{\log}(\rd)=\mathfrak{D}(\Phi_b)$ and thus $\mu\notin \mathfrak{D}(\Phi_b^{m+1})$. 
Suppose $n\geq 1$. Then $\Phi_b^n(\mu)$ is definable and equal to $\widetilde{\Phi}_{b,n}(\mu)$ by the assumption of induction. 
Denoting the L\'evy measure of $\Phi_b^n(\mu)$ by $\nu_n$, we have
\begin{align*}
\int_{|x|>1}\log_b|x|&\nu_n(dx)=\int_0^\infty dt\int_\rd \log_b^+ \left|b^{-[f_{n-1}^*(t)]}x\right|\nu(dx)\\
&=\int_{|x|>1}\nu(dx)\int_0^{f_{n-1}([\log_b|x|]+1)} \left(\log_b|x|-[f_{n-1}^*(t)]\right)dt\\
&\geq \int_{|x|>1}\nu(dx)\int_0^{(\log_b|x|)^n/n!} \left\{\log_b|x|-(n!\,t)^{\frac{1}{n}}\right\}dt\\
&=\int_{|x|>1}\frac{(\log_b|x|)^{n+1}}{(n+1)!}\nu(dx)=\infty.
\end{align*}
Thus, $\Phi_b^n(\mu)\notin I_{\log}(\rd)=\mathfrak{D}(\Phi_b)$ and hence $\Phi_b^{m+1}(\mu)$ is not definable. 
Therefore $\mathfrak{D}(\Phi_b^{m+1})\subset I_{\log^{m+1}}(\rd)$.
\end{proof}
Recall that the definition of nested subclasses of $L(b,\rd)$ in 
\citet{MaejimaNaito1998} mentioned in Introduction. 
The following theorem shows that $L_m(b,\rd)$ can be realized as $\mathfrak R (\Phi_b^{m+1})$.
\begin{theorem}
Suppose that $b>1$ and $m\in\Z_+$. Then,
$$
\Phi_b^{m+1}\left(I_{\log^{m+1}}(\rd)\right)=L_m(b,\rd).
$$
\end{theorem}
\begin{proof}
Let us show the statement by induction. If $m=0$, the assertion is Corollary \ref{s.s.d.}. 
Assume that the assertion is true for $m{-}1$ in place of $m$. 

We first show that $\Phi_b^{m+1}\left(I_{\log^{m+1}}(\rd)\right)\supset L_m(b,\rd)$. If $\mu\in L_m(b,\rd)$, 
there exists $\rho\in L_{m-1}(b,\rd)$ satisfying $\widehat{\mu}(z)=\widehat{\mu}(b^{-1}z)\widehat{\rho}(z)$. 
The assumption of induction implies that $\rho=\Phi_b^m(\rho_0)$ for some 
$\rho_0\in\mathfrak{D}(\Phi_b^m)$.
On the other hand, Theorem \ref{equivalence_between_decomposability_and_mapping} yields that
$\rho\in \mathfrak{D}(\Phi_b)$ and $\mu=\Phi_b(\rho)$.
Hence $\rho_0\in \mathfrak{D}(\Phi_b^{m+1})$ and $\mu=\Phi_b^{m+1}(\rho_0)$.
Thus $\mu\in\Phi_b^{m+1}\left(I_{\log^{m+1}}(\rd)\right)$.

To show the converse inclusion of two sets, suppose $\mu\in\Phi_b^{m+1}\left(I_{\log^{m+1}}(\rd)\right)$.
Then, $\mu=\Phi_b^{m+1}(\rho)=\Phi_b\left(\Phi_b^{m}(\rho)\right)$ for some $\rho\in\mathfrak{D}(\Phi_b^{m+1})$.
The assumption of induction implies that $\Phi_b^{m}(\rho)\in L_{m-1}(b,\rd)$, and Theorem 
\ref{equivalence_between_decomposability_and_mapping} yields that
$\widehat{\mu}(z)=\widehat{\mu}(b^{-1}z)\widehat{\Phi_b^{m}(\rho)}(z)$.
Thus $\mu\in L_m(b,\rd)$.
\end{proof}
Let $L_{\infty}(b,\rd)=\bigcap_{m=0}^\infty L_m(b,\rd)$. 
In \citet{MaejimaSatoWatanabe2000}, they studied $L_{\infty}(b,\rd)$ in the more general setting of 
operator semi-selfdecomposable distributions and as a special case, they proved that
$
L_{\infty}(b,\rd)=\overline{\mathit{SS}(b,\rd)},
$
where $\mathit{SS}(b,\rd)$ is the class of all semi-stable distributions with span $b$ and 
$\overline{\mathit{SS}(b,\rd)}$ 
denotes the closure of $\mathit{SS}(b,\rd)$ taken under convolution and weak convergence. 
What we want to emphasize here is that we have characterized $L_m(b,\rd)$ as the range of the mapping $\Phi_b^{m+1}$, 
and so we can conclude the following. 
Note that since $L_m(b,\rd)\supset L_{m+1}(b,\rd)$, $\lim_{m\rightarrow \infty }L_m(b,\rd)=L_\infty(b,\rd)$.
\begin{corollary}\label{closure_of_ss}
$$
\lim_{m\rightarrow \infty }\Phi_b^{m+1}\left(I_{\log^{m+1}}(\rd)\right)=L_{\infty}(b,\rd)=
\overline{\mathit{SS}(b,\rd)}.
$$
\end{corollary}
\begin{remark}
In \citet{MaejimaSato2009}, they proved that the limits of nested classes of several classes in 
$I(\rd)$ are identical 
with $L_{\infty}(\rd)$, which is known to be the same as the closure of the class of all stable 
distributions on $\rd$, 
$\overline{S(\rd)}$, say. 
Then a natural question arose.
Can we find mappings by which, as the limit of iteration, we get a larger or a smaller class than 
$\overline{S(\rd)}$? 
It is easy to see that $L_{\infty}(b,\rd)\supsetneqq L_{\infty}(\rd)$ so that $\overline{\mathit{SS}(b,\rd)}
\supsetneqq \overline{S(\rd)}$. 
\citet{Sato20072008} constructed mappings producing a class smaller than 
$\overline{S(\rd)}$.
Corollary \ref{closure_of_ss} shows that a mapping $\Phi_b$ produces a larger class than $\overline{S(\rd)}$ 
by iteration as a limit.
\end{remark}

\par\bigskip\noindent
{\bf Acknowledgment.}
The authors would like to thank Ken-iti Sato for his valuable comments.


\bibliographystyle{amsplain}

\end{document}